\documentclass[preprint,12pt]{elsarticle}

\usepackage{graphics}
\usepackage{graphicx}
\usepackage{epsfig}
\usepackage{anysize}
\usepackage{array}
\usepackage[fleqn]{amsmath}
\usepackage{latexsym,amssymb,amscd,amsthm,amsxtra}
\usepackage{subfigure}
\usepackage{tabularx}
\usepackage{enumerate}
\usepackage{multirow}
\usepackage{float}
\usepackage{placeins}
\usepackage{booktabs}
\usepackage{times}
\usepackage{soul,xcolor}
\usepackage{algorithm}
\usepackage{algpseudocode}
\usepackage[permil]{overpic}
\usepackage{tikz}
\usepackage{siunitx}
\usetikzlibrary{shapes.geometric, arrows}
\tikzstyle{startstop} = [rectangle, rounded corners, minimum width=2.5cm, minimum height=0.8cm, text centered, draw=black, fill=red!30, font=\small]
\tikzstyle{process} = [rectangle, minimum width=3.5cm, minimum height=0.8cm, text centered, draw=black, fill=blue!30, font=\small, text width=3.5cm, align=center]
\tikzstyle{arrow} = [thick,->,>=stealth]

\newtheorem{lemma}{Lemma}[section]
\newtheorem{definition}{Definition}[section]
\usepackage[most]{tcolorbox}

\newcommand{\blue}[1]{{\color{blue} #1}}
\newcommand{\green}[1]{{\color{green} #1}}
\newcommand{\magenta}[1]{{\color{magenta} #1}}

\usepackage{arydshln}

\journal{Engineering Analysis with Boundary Elements}

\def\R{{\mathbb{R}}}
\def\n{{\hat{n}}}
\begin{document}

\begin{frontmatter}

 \title{Efficient manifold evolution algorithm using adaptive B-Spline interpolation
 }
 \author{Muhammad Ammad\corref{cor1}\fnref{label1}}
\ead{21481199@life.hkbu.edu.hk}
\cortext[cor1]{Corresponding author}

\author{Leevan Ling\fnref{label1}}
\ead{lling@hkbu.edu.hk}

\affiliation[label1]{organization={Department of Mathematics, Hong Kong Baptist University},
             country={Hong Kong}}

\begin{abstract}
This paper explores an efficient Lagrangian approach for evolving point cloud data on smooth manifolds. In this preliminary study, we focus on analyzing plane curves, and our ultimate goal is to provide an alternative to the conventional radial basis function (RBF) approach for manifolds in higher dimensions. In particular, we use the B-Spline as the basis function for all local interpolations. Just like RBF and other smooth basis functions, B-Splines enable the approximation of geometric features such as normal vectors and curvature. Once properly set up, the advantage of using B-Splines is that their coefficients carry geometric meanings. This allows the coefficients to be manipulated like points, facilitates rapid updates of the interpolant, and eliminates the need for frequent re-interpolation. Consequently, the removal and insertion of point cloud data become seamless processes, particularly advantageous in regions experiencing significant fluctuations in point density. The numerical results demonstrate the convergence of geometric quantities and the effectiveness of our approach. Finally, we show simulations of curvature flows whose speeds depend on the solutions of coupled reaction-diffusion systems for pattern formation.

\bigskip
\noindent \textbf{Research highlights}:
\begin{itemize}
    \item We employ B-Spline as the basis functions for local interpolations, enabling rapid updates of the interpolant without the need for re-interpolation at each time step.
    \item The use of B-Spline allows us to assign geometric meanings to control points, facilitating the movement of these coefficients as if they were points.
    \item Our approach provides precision in estimating normal vectors and curvature vectors that is comparable to the conventional PHS-RBF+Poly method.
    \item We have introduced adaptive mechanisms for adding and removing points from the point cloud in response to significant changes in point density.
    \item Numerical results validate the effectiveness of our method in simulating curvature-driven flows and handling complex geometries in coupled reaction-diffusion systems.
\end{itemize}

\end{abstract}


\end{frontmatter}

\section{Introduction}\label{sec:intro}
In this study, we develop a computational framework to model the evolution of manifolds using a controlled evolution mechanism. By defining a speed function normal to the evolving manifold, which depends on intrinsic properties such as curvature and is influenced by the solutions of underlying partial differential equations (PDEs), we can accurately simulate the dynamic changes of domains over time. This approach is crucial to understanding and predicting complex behaviors in physical and biological systems, where domain evolution plays a fundamental role in underlying dynamics \cite{yang2006embedded,maitre2009applications, eilks2008numerical, hegarty2021modelling,venkataraman2012adaptive,fritz2023tumor}.

The study of manifold evolution, particularly in higher dimensions, has a rich history in geometric modeling and computational geometry. Foundational work has demonstrated how geometric shapes dynamically evolve under the influence of intrinsic properties, such as curvature or external forces, while maintaining essential geometric constraints such as the preservation of volume, area, or other characteristics \cite{gage1986area, ma2014non,tsai2015length}. For example, Gage and Hamilton's seminal work on the evolution of convex curves \cite{gage1986heat} revealed how these curves evolve, inspiring a wide range of studies on the evolution of more complex manifolds under various geometric and physical constraints. Since then, these approaches have been extended to non-convex and immersed manifolds, where local and global properties influence the evolution of the manifold \cite{sesum2020evolution, wang2023evolution, wang2018evolution,dittberner2021curve}.

A particular area of interest within manifold evolution is the study of planar or space curves and their dynamics under curvature-driven flows. For instance, mean curvature flows generate minimal surfaces from curves with non-vanishing torsion, highlighting the significant roles that curvature and torsion play in their evolution \cite{minarcik2022minimal}. The dynamics of specific curves, such as helix curves, can be analyzed through their velocity and geometric properties \cite{abdel2014evolution}. However, as these curves evolve, one of the major challenges lies in maintaining control over the shape, particularly in preserving intricate geometric details. To address this, research has shown that introducing additional constraints during the evolution process can effectively prevent distortions or instabilities, ensuring the preservation of desired features \cite{nagasawa2019interpolation, gao2023star}. These constraints are essential for managing more complex shapes and are highly relevant to many real-world applications \cite{jiang2008non}, such as fluid dynamics and biology, where evolving filaments and membranes respond to curvature or external forces \cite{guan2009quermassintegral, huisken2001inverse}.

Recently, advanced computational methods, including the level set method (Eulerian approach), have been used to model the evolution of 1D and 2D manifolds, particularly in scenarios involving topology changes, such as splitting or merging \cite{kim2020pattern, yu2020survey}. However, our approach adopts a Lagrangian framework, which allows for direct tracking of individual manifold points and provides higher accuracy in capturing geometric deformations. This direct tracking makes the Lagrangian approach particularly advantageous for problems where pointwise details and local geometric properties are critical to the analysis.

To organize our contributions, we first present the general Lagrangian framework for moving unknown manifolds in Section \ref{sec:general_framework}, where we describe the evolution of a point cloud representing a smooth manifold. The motion is driven by a velocity field, typically influenced by geometric properties such as curvature and the normal vector. Section \ref{sec:adaptive_bspline} introduces our adaptive method for point cloud evolution using B-Spline basis functions. We detail our approach to B-Spline curve interpolation and control point refinement, highlighting its flexibility in approximating evolving geometries. In Section \ref{sec:results}, we present results and discussion, examining the impact of stencil size, B-Spline degree, and point cloud density on the interpolation accuracy. We conclude with a simulation of interface dynamics, including an evolving circle boundary coupled with reaction-diffusion processes. Finally, in Section \ref{sec:conclusion}, we summarize our findings and propose directions for future research.
\section{Lagrangian Framework for Moving Curves} \label{sec:general_framework}
Let \(\mathcal{M}(t)\) represent a smooth one-dimensional evolving curve embedded in \(\mathbb{R}^d\) over the time interval \( t \in [0, T] \). The motion of any point \(\vec{q}(t) \in \mathcal{M}(t)\) on this curve can be described by
\begin{equation}\label{dq/dt}
  \frac{d\vec{q}}{dt} = \vec{V}(\vec{q},t),
\end{equation}
where \(\vec{V}(\vec{q},t)\) is a predefined velocity field.
We consider the discrete cases involving the evolution of a point cloud
\begin{equation}\label{Q(t)}
  Q(t) := \{\vec{q}_i(t)\}_{i=1}^{N} \subset \mathcal{M}(t)
\end{equation}
on the initial curve with $t=0$, often referred to as an ``unknown manifold'' \cite{yan2023spectral,liang2024solving} in recent literature.
Although we may have analytic information about the initial curve \(\mathcal{M}(0)\), this information can be used to compute geometric quantities at the points in \(Q(0)\). However, as the curve evolves to $t>0$, the shape of \(\mathcal{M}(t)\) becomes unknown, and therefore the geometric properties at the points in \(Q(t)\) must be approximated numerically.

Consider the case of standard \emph{curvature motion} of a curve, characterized by the velocity field \(\vec{V}(\vec{q}, t) = -\kappa(\vec{q}) \mathbf{n}(\vec{q})\), where \(\kappa(\vec{q})\) represents the curvature at point \(\vec{q}(t)\), and \(\mathbf{n}(\vec{q})\) denotes the outward normal vector at the same point.
Under this curvature-driven flow, the curve tends to shrink, with regions of higher curvature shrinking more rapidly.

This phenomenon is often observed in the process of surface regularization and simplification \cite{petras2019least}.
To update the point cloud from \(Q(t)\) to \(Q(t + \Delta t)\), we will apply the Lagrangian framework by discretizing \eqref{dq/dt} to move each data point. Suppose the velocity field \(\vec{V}(\vec{q}_i, t)\) is known for each point \(\vec{q}_i(t)\):
\[
\vec{V}_i(\vec{q}, t) = -\kappa_i(t) \n_i(t),
\]
where \(\kappa_i(t)\) is the curvature, and \(\n_i(t)\) is the normal vector at \(\vec{q}_i(t)\). Estimating these quantities is the focus of this paper and will be discussed in the following. For now, we assume that the velocity field \(\vec{V}\) can be evaluated wherever it is needed. The position of each point is then updated by discretizing the velocity field in time. The new position at the next time step \(t + \Delta t\) is given by:
\begin{equation}\label{eq:point_evolution}
    \vec{q}_i(t + \Delta t) = \vec{q}_i(t) + \Delta t \cdot \vec{V}_i(\vec{q}, t).
\end{equation}
The formulation provided offers a straightforward method for updating point positions; however, practical challenges arise when estimating geometric quantities, adapting the representation of the evolving curve to account for local variations in geometry, and ensuring computational efficiency. To address these challenges, we introduce an adaptive method based on B-Spline interpolation, which dynamically refines the point cloud representation to ensure accuracy and smoothness. A detailed description of this adaptive method, along with its implementation, is provided in the next section.

\section{Adaptive method for point cloud evolution using B-Splines} \label{sec:adaptive_bspline}
The study of dynamic interfaces and evolving geometries requires computational tools capable of real-time adaptation to changes in model parameters and structures. In this work, we present a new approach for evolving point cloud data using adaptive B-Spline interpolation, which ensures accurate curve representations by adjusting control points and preserving smoothness, even in regions with complex geometric variations.

A key feature of our framework is that local B-spline interpolation is not only used to fit the evolving point cloud and thus serves as the basis for directly and efficiently computing geometric quantities (such as normals and curvature) at all core points, but also provides additional computational tools for managing the point cloud over time. Unlike traditional approaches that require re-interpolating or recalculating all coefficients at each time step, our method evolves the B-Spline control points (i.e., interpolation coefficients) alongside the data points, providing a better initial guess for iterative methods. Furthermore, the local refinement capability of B-Spline interpolation enables the efficient insertion and removal of points, allowing the representation to adapt dynamically to changes in point density or geometric complexity.

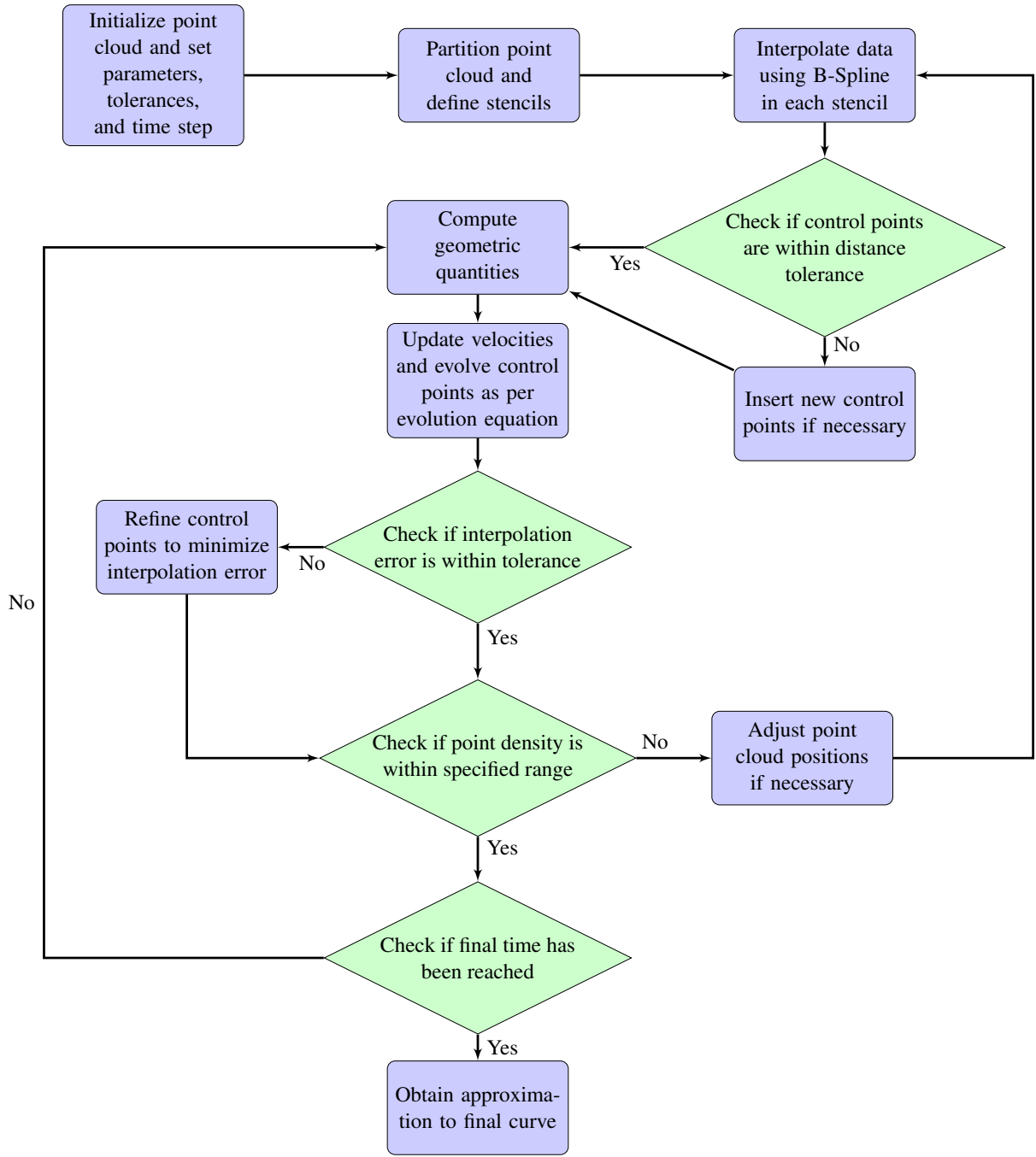
\begin{figure}[pt!]
\centering
\resizebox{\textwidth}{!}{
\begin{tikzpicture}[node distance=3cm, auto, every node/.style={font=\small}]
\tikzstyle{block} = [rectangle, draw, fill=blue!20, text width=7em, text centered, rounded corners, minimum height=4em]
\tikzstyle{decision} = [diamond, aspect=2, draw, fill=green!20, text width=9em, text badly centered, inner sep=1pt]
\tikzstyle{line} = [draw, very thick, -latex']

\node [block] (init)
{Initialize point cloud and set parameters, tolerances, and time step};

\node [block, right of=init, node distance=6cm] (partition)
{Partition point cloud and define stencils};

\node [block, right of=partition, node distance=6cm] (fit_curve)
{Interpolate data using B-Spline in each stencil};

\node [decision, below of=fit_curve, node distance=3.1cm] (check_distance)
{Check if control points are within distance tolerance};

\node [block, below of=check_distance, node distance=3cm] (insert_control_points)
{Insert new control points if necessary};

\node [block, left of=check_distance, node distance=6.2cm] (calculate_geom)
{Compute geometric quantities};

\node [block, below of=calculate_geom, node distance=2.4cm] (update_points)
{Update velocities and evolve control points as per evolution equation};

\node [decision, below of=update_points, node distance=3.cm] (check_error)
{Check if interpolation error is within tolerance};

\node [block, left of=check_error, node distance=5.2cm] (optimize_control_points)
{Refine control points to minimize interpolation error};

\node [decision, below of=check_error, node distance=3.8cm] (check_dist_qi)
{Check if point density is within specified range};

\node [block, right of=check_dist_qi, node distance=5.8cm] (redistribute_qi)
{Adjust point cloud positions if necessary};

\node [decision, below of=check_dist_qi, node distance=3.6cm] (converged)
{Check if final time has been reached};

\node [block, below of=converged, node distance=2.7cm] (end)
{Obtain approximation to final curve};

\path [line] (init) -- (partition);
\path [line] (partition) -- (fit_curve);
\path [line] (fit_curve) -- (check_distance);

\path [line] (check_distance) -- node [near start] {Yes} (calculate_geom);
\path [line] (check_distance) -- node [near start, right] {No} (insert_control_points);
\path [line] (insert_control_points) -- (calculate_geom);

\path [line] (calculate_geom) -- (update_points);
\path [line] (update_points) -- (check_error);

\path [line] (check_error) -- node [near start, right] {Yes} (check_dist_qi);
\path [line] (check_error) -- node [near start, below] {No} (optimize_control_points);
\path [line] (optimize_control_points) |- (check_dist_qi);

\path [line] (check_dist_qi) -- node [near start, above] {No} (redistribute_qi);
\path [line] (check_dist_qi) -- node [near start, right] {Yes} (converged);

\path [line] (redistribute_qi.east) --++(2.5,0)|- (fit_curve);
\path [line] (converged) -- node {Yes} (end);
\path [line] (converged.west) -- ++(-5,0) |- node [near start, left] {No} (calculate_geom);

\end{tikzpicture}}
\caption{This flowchart provides a detailed workflow for the evolution of a point cloud. It includes steps for initialization, partitioning, interpolation using B-Splines, checking tolerances, updating point positions, and refining control points, culminating in obtaining an approximation to the final curve.}\label{fig:flowchart}
\end{figure}

Our approach employs a partitioned and overlapping stencil strategy to enhance efficiency and precision while maintaining critical geometric properties. By iteratively refining the point cloud, the method achieves accurate representation and smooth transitions, addressing key challenges in point cloud evolution and providing a foundation for future applications across scientific and engineering fields.
For a clearer understanding of the discussed concepts, a visual summary of the algorithm is presented in Figure~\ref{fig:flowchart}. Additionally, a comprehensive terminology map comparing B-Splines and RBFs is provided in Table~\ref{tab:bspline_rbf}.

\begin{table}[t]
\centering
\setlength{\tabcolsep}{18pt} 
\begin{tabular}{@{}ll@{}}
\toprule
\textbf{B-Spline Terminology}       & \textbf{RBF Terminology}                     \\ \midrule
Knots                               & Centers of RBF                               \\
Control points                      & Coefficients in RBF Representation           \\
Degree of Spline                    & Order of the RBF                             \\
Spline Space                        & Native Space of the RBF                      \\
Spline Interpolation                & Kernel Interpolation                         \\
\midrule
\multicolumn{2}{c}{\textbf{Hyperparameters}} \\
\midrule
Degree, Knots                       & Kernel type, Shape parameter        \\
\bottomrule
\end{tabular}
\caption{Comparison of terminologies used in B-Splines and Radial Basis Functions (RBFs), highlighting parallels in concepts and applications.}\label{tab:bspline_rbf}
\end{table}

\subsection{Designing a cover for point clouds}\label{sec:stencil}
Assume that we have an initial point cloud \( Q(0)\), as defined in \eqref{Q(t)}, sampled from an initial curve \( \mathcal{M}(0) \subset \mathbb{R}^2 \). To begin the computation of local geometric quantities, \( Q(0) \) is partitioned into distinct, connected but non-overlapping \emph{core covers}
 \(\{ \mathcal{N}_{k}^\text{core} \}_k \), each containing \( m_{c,k} \) points and $\sum_k m_{c,k} = N$.
To ensure smooth transitions and stable computations near boundaries, each core cover is expanded by incorporating
\( m_{b,k} \) points from adjacent areas, creating a \emph{boundary subset}
\(\{  \mathcal{N}_{k}^\text{bdy} \}_k\).
The combined structure gives the \emph{stencil} for local interpolation and is described by:
\begin{equation}
    \mathcal{N}_{k} = \mathcal{N}_{k}^\text{core} \cup \mathcal{N}_{k}^\text{bdy},
    \label{eq:cover}
\end{equation}
where \( m_k=m_{c,k}+m_{b,k} \) represents the total number of points in a stencil, and \( \mathcal{N}_{m_{c,k}}^\text{core} \cap \mathcal{N}_{m_{b,k}}^\text{bdy} = \emptyset \).
\[
Q(t) = \bigcup_k \mathcal{N}_{k}^\text{core} \subseteq \bigcup_k \mathcal{N}_{k}.
\]
This construction ensures that every point in $\vec q_i(t) \in Q(t)$  is included in a unique core cover
label $k=k(i)$ and, hence,a unique stencil, denoted by $\mathcal{N}_{k}$, for the approximation of the function at $\vec q_i(t)$.

\begin{figure}[t!]
\centering
        \begin{overpic}[width=10cm]{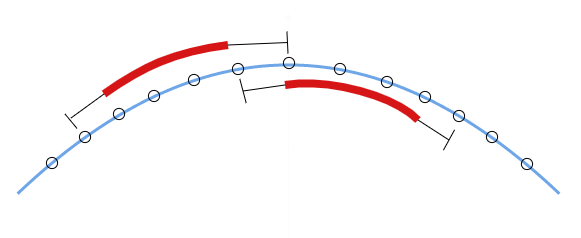}
              \put(210,310)  {\rotatebox{30}{\scriptsize $\mathcal{N}_{k}^\text{core}$}}
              \put(560,230)  {\rotatebox{340}{\scriptsize $\mathcal{N}_{k+1}^\text{core}$}}
        \end{overpic}
\caption{Partitioning of point cloud data into non-overlapping cores (red segments), which consist of subsets of consecutive points, and their associated stencils (black lines), which extend beyond the cores by including additional boundary points to ensure a smooth transition across the curves.}
\label{fig:cover3.1}
\end{figure}

The visualization in Figure~\ref{fig:cover3.1} demonstrates how the point cloud is partitioned into distinct regions. Each core cover consists of a subset of consecutive points, while the stencils incorporate additional boundary points on either side of the core covers to ensure continuity. This approach allows the stencils to cover a larger neighborhood of the curve for analysis.

The choice of stencil size \( m_k \) balances computational complexity and accuracy; a smaller \( m_k \) reduces complexity but may compromise accuracy, while a larger \( m_k \) increases accuracy at the cost of higher computational demand. The designed stencils and their partitions are crucial for estimating geometric quantities and ensuring smooth transitions across overlapping neighborhoods.

\subsection{Interpolation of parametric curves}
 \label{subsec:spline_interpolation}
With the stencil defined, we now have the domains to interpolate the underlying geometry. The next step involves parameterizing the data within each cover and constructing interpolated curves. This section focuses on constructing parametric curves for each cover independently using B-Splines, with geometric significance assigned to the control points. For simplicity, we restrict our explanation to a single extended cover, as the same methodology applies independently to all covers.

Although the interpolation concept aligns with simpler methods, it involves treating each coordinate independently. We consider the local stencil introduced in \eqref{eq:cover} 
and drop the subscript~$k$ for stencil label, as this does not cause any confusion.
The component-wise interpolation data are given by:

\[
\big\{ \big( u_i, [ \vec f(u_i) ]_j \big) \big\}_{i=1}^m :=
\big\{ \big( u_i, [ \vec q_i ]_j \big) \big\}_{i=1}^m \subset \R \times \R
\quad \text{for } j=1,2,
\]
where data site \( u_i \) are parameter values to be assigned to the data points \( \vec{q}_i = (q_{i,1}, q_{i,2})^T \), and  $\vec{f}(u) = [f_1(u), f_2(u)]^T$ is the unknown curve to be interpolated. To assign parameter values \( u_i \) to the data points, we employ the traditional chord length parameterization method \cite{piegl2012nurbs}.
Specifically, the parameterization begins by setting \( u_1 = 0 \) for the starting point and \( u_n = 1 \) for the ending point. The
parameter values \( u_i \), for \( i = 2, \dots, m-1 \) are determined on the basis of the cumulative chord lengths between consecutive points as follows:

\[
d = \sum_{i=2}^{m} d_i := \sum_{i=2}^{m} \lVert \vec{q}_i - \vec{q}_{i-1} \rVert_2,
\]
Namely, we set
$u_i = u_{i-1} + {d_i}/{d} \in (0,1)$ for $i = 2, \dots, m-1$.

\newcommand\solidrule[1][1cm]{\rule[0.5ex]{#1}{.4pt}}
\newcommand\dashedrule{\mbox{%
  \solidrule[2mm]\hspace{2mm}\solidrule[2mm]\hspace{2mm}\solidrule[2mm]}}
\begin{figure}[t!]
\centering
        \begin{overpic}[width=9cm]{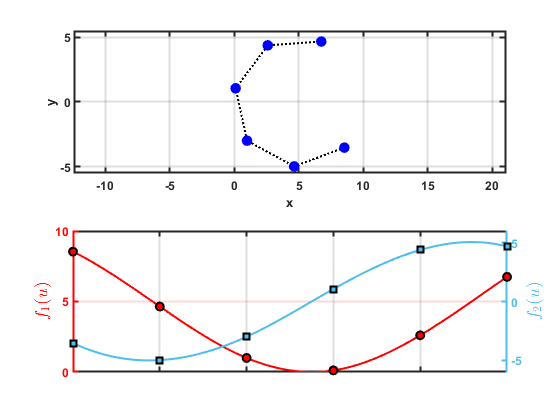}
              \put(560,460) {\scriptsize \magenta{$d_1$}}
              \put(465,465) {\scriptsize \magenta{$d_2$}}
              \put(415,530) {\scriptsize \magenta{$d_3$}}
              \put(435,625) {\scriptsize \magenta{$d_4$}}
              \put(515,660) {\scriptsize \magenta{$d_5$}}
            \put(600,505) {\scriptsize $q_1$}
             \put(510,470) {\scriptsize $q_2$}
             \put(450,500) {\scriptsize $q_3$}
             \put(435,585) {\scriptsize $q_4$}
             \put(480,640) {\scriptsize $q_5$}
             \put(570,650) {\scriptsize $q_6$}
            \put(120,50) {\scriptsize $u_1$}
            \put(270,50) {\scriptsize \rotatebox{0}{$u_2$}}
            \put(425,50) {\scriptsize \rotatebox{0}{$u_3$}}
            \put(580,50) {\scriptsize \rotatebox{0}{$u_4$}}
            \put(735,50) {\scriptsize \rotatebox{0}{$u_5$}}
            \put(880,50) {\scriptsize \rotatebox{0}{$u_6$}}
        \end{overpic}
\caption{Illustration of parameterization and interpolation in the given data. (Top) The input data points, denoted as $q_i$ and labeled as \blue{$\bullet$}, with $d_i$ representing the distances between consecutive points. (Bottom) The parameterization values, $u_1, u_2, \dots, u_6$, are computed based on chord length. The red curve represents the smooth interpolation of the x-coordinates, with discrete points labeled as \textbf{$\circ$}, while the blue curve represents the smooth interpolation of the y-coordinates, with discrete points labeled as \textbf{$\square$}.}
\label{fig:interpol_curve}
\end{figure}
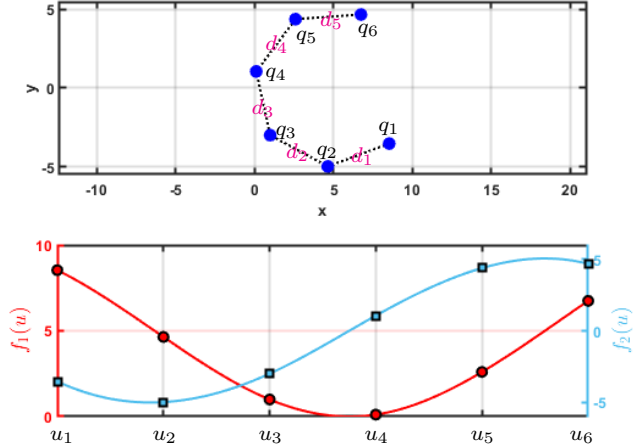

This approach of picking data points normalizes the parameter values to the interval \([0, 1]\) and ensures that they are distributed proportionally to the distances between consecutive points along the data. Figure \ref{fig:interpol_curve} illustrates this parameterization process. The top panel displays the input data points and their geometric relationships, highlighting the distances between consecutive points. The bottom panel visualizes the parameterization values and their role in capturing the structure of the data, enabling smooth interpolation of the x- and y-coordinate components independently.

Following parameterization, the next step is to construct the interpolant.  If we adopt an identical ansatz for both interpolants \( f_j: \R \to \R \) of the curve parameterization vector function \( \vec{f} := [f_1, f_2]^T \), then we will end up with a coefficient vector \( \vec{P}_\ell \) for each basis function \( \phi_\ell \). We numerically expand each component function, for \( j = 1, 2 \), in the form
\[
f_j(u) = \sum_{\ell=1}^m [\vec{P}_\ell]_j \phi_\ell(u),
\]
and impose the interpolation conditions at the parameter values \( \{u_i\}_{i=1}^m \) to yield
\[
f_j(u_i) = \sum_{\ell=1}^m [\vec{P}_\ell]_j \phi_\ell(u_i) = [\vec{q}_i]_j,
\quad \text{for } i = 1, \ldots, m.
\]
Solving these two interpolation problems of size proportional to the number of interpolation points \( m \) allows us to determine the set of control point vectors: $
\vec{P}_\ell = \big( [\vec{P}_\ell]_1, [\vec{P}_\ell]_2 \big)^T \in \R^2$, and express the vector interpolant function in the following vector form:

\[
\vec{f}(u) = \sum_{\ell=1}^m \vec{P}_\ell \phi_\ell(u) : \R \to \R^2.
\]
When the basis functions \( \phi_\ell \) are chosen to be B-spline basis functions, the interpolation problem becomes B-spline interpolation. In this case,  the B-Spline basis function $ \phi_j^p: \mathbb{R} \to \mathbb{R} $, for $ m $ control points of degree $ p $, over a knot vector $ \Psi = [\psi_1, \dots, \psi_{m+p+1}]^T \in \mathbb{R}^{m+p+1} $, is recursively defined by the Cox-de Boor formula:
\begin{equation*}
\phi_i^0(u) =
\begin{cases}
1, & \text{if } u \in [\psi_j, \psi_{j+1}) \\
0, & \text{otherwise}
\end{cases}
\end{equation*}

\begin{equation*}
\phi_j^p(u) = \frac{u - \psi_j}{\psi_{j+p} - \psi_i} \phi_i^{p-1}(u) + \frac{\psi_{j+p+1} - u}{\psi_{j+p+1} - \psi_{j+1}} \phi_{j+1}^{p-1}(u),
  \label{eq:basis_recursive}
\end{equation*}
for $ j = 1, \dots, n $ and $ p \geq 1 $.

The basis functions $ \phi_j^p(u) $ in B-Spline interpolation determine how each control point $ \vec{P}_j $ influences the curve at a given parameter value $ u $. These functions distribute the control points' effects across the (parameter) domain. Two common types of B-Spline basis functions are \emph{periodic} and \emph{open uniform}.

\begin{figure}
    \begin{center}
        \subfigure[]{\includegraphics[width=2.8in]{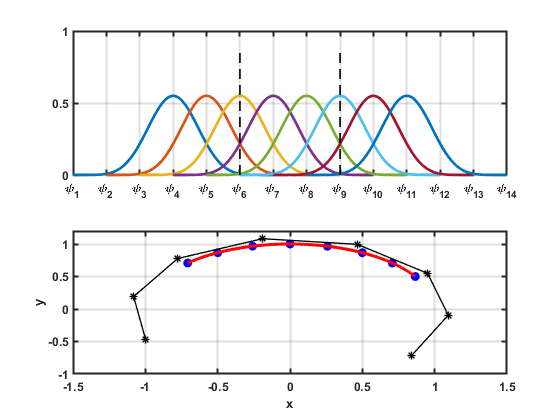}}
        \subfigure[]{\includegraphics[width=2.8in]{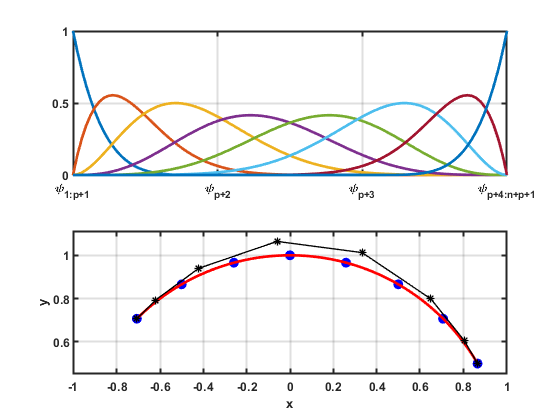}}
    \caption{Comparison of B-Spline basis functions and their relation to interpolation (a) Periodic B-Spline basis functions (top) provide smooth cyclic transitions across the knot vector domain. The associated curve (bottom) passes through the extended cover points, but the first and last control points lie outside the data domain, complicating the computation of geometric quantities at the boundaries. (b) Open uniform B-Spline basis functions (top) result in a curve (bottom) that passes through all extended cover points, with control points lying entirely within the data domain, including boundary points, facilitating geometric computations.}
    \label{fig:basisfunction}
    \end{center}
\end{figure}

\begin{definition}[\cite{marsh2005applied}]\label{def:periodicbasis}
A periodic B-spline of degree \( p \) and with a control point \( n \) is obtained by arbitrarily choosing the knots \( \psi_ \leq \cdots \leq \psi_n \) and then setting
\[\psi_{m+j} = \psi_{m+j-1} + (\psi_j - \psi_{j-1}),\]
for \( j = 1, \dots, p+1 \). A knot vector of this form is called a periodic knot vector. In particular, a uniform B-spline is a special case of a periodic B-spline.
\end{definition}

\begin{definition}[\cite{marsh2005applied}]\label{def:openuniform}
An open uniform B-Spline is a B-Spline where the knot vector is uniform except at its two ends, where the knot values are repeated \( k \) times, with \( k = p + 1 \). For a given \( m \) (number of control points) and \( k \), the center vector can be defined as:

\[
\psi_j =
\begin{cases}
0, & 0 \leq j < k, \\
j - k + 1, & k \leq j \leq m, \\
m - k + 2, & m < j \leq m + k.
\end{cases}
\]
This construction ensures that the B-Spline curve interpolates the first and last control points, providing precise control at the endpoints.
\end{definition}
We remark that the selection of B-Spline basis functions—whether periodic or open-uniform depends on the specific application requirements. Periodic B-Spline basis functions, as shown in Figure \ref{fig:basisfunction}(a), provide cyclic continuity by repeating themselves across the knot vector. This inherent property ensures smooth and seamless transitions across the entire domain, making them particularly suitable for applications involving closed manifolds or periodic structures such as loops, surfaces, or cyclic phenomena. However, their lack of precise endpoint control, as they do not satisfy \( \vec{f}(0) = \vec{P}_{1} \) and \( \vec{f}(1) = \vec{P}_{m} \), limits their use in scenarios where accurate boundary conditions are required.

In contrast, open uniform B-Spline basis functions, as demonstrated in Figure \ref{fig:basisfunction}(b), are ideal for applications requiring precise control over the curve endpoints. These basis functions guarantee that the curve interpolates the first and last control points, satisfying \( \vec{f}(0) = \vec{P}_{1} \) and \( \vec{f}(1) = \vec{P}_{m} \). This property makes them particularly effective for tasks such as dynamic point cloud evolution, where accurate interpolation and boundary control are vital to preserving geometric fidelity and maintaining stability.

In this work, we adopt open uniform B-Spline basis functions due to their ability to provide local control and endpoint accuracy, which are essential for evolving point clouds and ensuring stable boundary conditions. However, in future work, periodic B-Spline basis functions will be considered for applications requiring global interpolation of closed manifolds, where smooth cyclic continuity is necessary to handle the geometry of closed shapes effectively.

\subsection{Distance-based Control Point Refinement}
\label{subsec:control_point_refinement}
In the interpolation process, the given point cloud data \( Q(t) \) produces a smooth curve for each stencil, with the control points \( \{\vec{P}_j\}_{j=1}^m \) fixed as a result. However, to ensure that the control points accurately capture the geometry of the curve, it is necessary to measure how closely the curve aligns with its control points. If the control points are too far from the curve, they may fail to adequately represent the curve's geometry, leading to inaccuracies in subsequent geometry-driven updates. To address this, we refine the control points to better reflect the geometric properties of the curve.

To formalize this refinement process, we adopt the following definition for measuring the distance between the curve and its control points, as proposed in \cite{scharf2003computing}:

\begin{definition}[\cite{scharf2003computing}]
Given a control point \( \vec{P}_j = (x_j, y_j) \) and a degree-\(p\) B-Spline curve \( \vec{f}(u) = (x(u), y(u)) \) for \( u \in [0, 1] \), let \( u_1, u_2, \dots, u_m \) denote the odd-multiplicity roots of the polynomial:
\[
P(u) = (x_j - x(u)) \cdot x'(u) + (y_j - y(u)) \cdot y'(u),
\]
which is of degree \( 2p - 1 \). Let \( u_0 = 0 \) and \( u_{m+1} = 1 \). The distance between \( \vec{f}(u) \) and \( \vec{P}_j \) is then computed as:
\[
\text{dist}_f(\vec{P}_j) = \min_{0 \leq k \leq m+1} \| \vec{f}(u_k) - \vec{P}_j \|.
\]
\end{definition}
Using this distance, we introduce the deviation metric, which identifies the control point farthest from the curve:
\[
\epsilon(\{\vec{P}_j\}, \vec{f}) = \max_{1 \leq j \leq m} \text{dist}_f(\vec{P}_j).
\]
This metric quantifies the largest deviation between the curve and its control points, providing a measure of how accurately the control polygon represents the geometry of the curve. If the deviation metric \( \epsilon(\{\vec{P}_j\}, \vec{f}) \) exceeds a prescribed tolerance \( \epsilon_{\text{tol}} \), refinement of the control points is required to improve the representation.

Refinement is achieved through knot insertion, a standard technique that adjusts the control polygon while preserving the smoothness of the curve. The refinement process exploits the nesting property of spline spaces, as described in the following lemma:
\begin{lemma}[\cite{lyche2008spline}]
\label{lemma:knotinsertion}
Let \( p \) be a positive integer and let \( \Psi \) be a knot vector containing at least \( p + 2 \) knots. If \( \tilde{\Psi} \) is a knot vector such that \( \Psi \subseteq \tilde{\Psi} \), then the spline space associated with \( \Psi \) is a subspace of the spline space associated with \( \tilde{\Psi} \), i.e.,
$\mathbb{S}_{p, \Psi} \subseteq \mathbb{S}_{p, \tilde{\Psi}}.$
\end{lemma}

Using the property $\mathbb{S}_{p, \Psi} \subseteq \mathbb{S}_{p, \tilde{\Psi}}.$ of the lemma~\ref{lemma:knotinsertion}, where \( \mathbb{S}_{p, \Psi} = \text{span}\{\phi_{1,\Psi}^p, \dots, \phi_{n,\Psi}^p\} \) and \( \mathbb{S}_{p, \tilde{\Psi}} = \text{span}\{\phi_{1,\tilde{\Psi}}^p, \dots, \phi_{m,\tilde{\Psi}}^p\} \), we introduce a new knot \( \psi^* \) into the knot vector \( \Psi \) at the parameter value \( u^* \), where the maximum deviation $\epsilon_{\text{max}} = \max_{u \in [0, 1]} \| \vec{f}(u) - \vec{\Gamma}(u) \|$ is achieved. The resulting knot vector \( \tilde{\Psi} = \Psi \cup \{\psi^*\} \) expands the spline space, allowing the control polygon to better approximate the geometry of the curve.

\begin{figure}[t]
    \centering
    \includegraphics[width=3in]{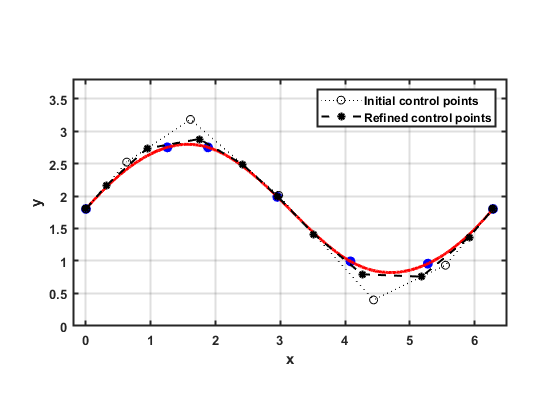}
    \caption{A demonstration of control point refinement using knot insertion. The original control points \( \vec{P}_j \) are labeled as \( \circ \), and the refined points \( \vec{P}_j^{\prime} \), obtained through knot insertion, are labeled as \( * \). This refinement results in a control polygon more closely aligned with the B-Spline curve.}
    \label{fig:knotinsertion}
\end{figure}

After refinement, the updated B-Spline curve is expressed in terms of the new knot vector \( \tilde{\Psi} \):
\[
\vec{f}(u) = \sum_{j=1}^m \vec{P}_j \phi_j^p(u) = \sum_{j=1}^{m+1} \vec{P}_j^{\prime} \phi_j^p(u),
\]
where the refined control points \( \vec{P}_j^{\prime} \) are calculated using the recurrence relation:
\[
\vec{P}_j^{\prime} = (1 - \alpha_j) \vec{P}_{j-1} + \alpha_j \vec{P}_j,
\]
and the weight \( \alpha_j \) is defined as:
\[
\alpha_j = \frac{\psi^* - \psi_j}{\psi_{j+p} - \psi_j}.
\]
This localized refinement ensures that the control polygon adapts to the geometry of the curve in regions where \( \epsilon(\{\vec{P}_j\}, \vec{f}) \) is large, reducing the deviation metric \( \epsilon(\{\vec{P}_j\}, \vec{f}) \) while preserving the smoothness of the curve. As shown in Figure~\ref{fig:knotinsertion}, the initial control points \( \vec{P}_j \), labeled by \( \circ \), are refined to updated control points \( \vec{P}_j^{\prime} \), labeled by \( * \), resulting in a control polygon that aligns more closely with the curve geometry and satisfies the prescribed tolerance \( \epsilon_{\text{tol}} \).

By iterative application of this process, the control polygon achieves a closer alignment with the B-Spline curve, reducing the deviation metric \( \epsilon(\{\vec{P}_j\}, \vec{f}) \) to within the given tolerance while maintaining the curve's smoothness.
\subsection{Curvature and normal computation}\label{sec:nandkappa}
After fitting the B-spline to the points in each stencil \( \mathcal{N}_{k} \) and refining it to meet the desired tolerance, we compute the geometric properties at each core point within the core subset \( \mathcal{N}_{k}^{\text{core}} \), specifically the normal vector \(\n(\vec{q}_i)\) and the curvature \(\kappa(\vec{q}_i)\).

For any point \(\vec{q}_i \in \mathcal{N}_{k(i)}^{\text{core}}\), the tangent vector at that point is derived from the first derivative of the B-spline curve \(\vec{f}(u)\). The tangent vector is expressed as:
\[
\vec{t}(u_i) = \vec{f}'(u_i) = \begin{bmatrix} f_1'(u_i) \\ f_2'(u_i) \end{bmatrix},
\]
The normal vector \(\n(\vec{q}_i)\) is orthogonal to \(\vec{t}(u_i)\), and provides the direction perpendicular to the curve at \(\vec{q}_i\). To compute it, we normalize the perpendicular direction of the tangent vector \(\vec{t}(u_i)\). The normal vector is given by:
\[
\n(\vec{q}_i) = \frac{1}{\|\vec{t}(u_i)\|} \begin{bmatrix} -f_2'(u_i) \\ f_1'(u_i) \end{bmatrix}.
\]
The curvature \(\kappa(\vec{q}_i)\) measures how sharply the curve bends at the point \(\vec{q}_i\). It is defined by:
\[
\kappa(\vec{q}_i) = \frac{\big|f_1'(u_i) f_2''(u_i) - f_2'(u_i) f_1''(u_i)\big|}{\big(f_1'(u_i)^2 + f_2'(u_i)^2\big)^{3/2}}.
\]
where,
\[
\vec{f}'(u) = \frac{d\vec{f}(u)}{du} = \begin{bmatrix} f_1'(u) \\ f_2'(u) \end{bmatrix},
\quad
\vec{f}''(u) = \frac{d^2\vec{f}(u)}{du^2} = \begin{bmatrix} f_1''(u) \\ f_2''(u) \end{bmatrix}.
\]

These calculations are carried out for all points \(\vec{q}_i \in \mathcal{N}_{k(i)}^{\text{core}}\) and are stored for use in the subsequent evolution of the point cloud.
\subsection{Geometry-driven point evolution and control point optimization}\label{sec:po_evol_intererror}
\begin{figure}
\centering
       \includegraphics[width=3in]{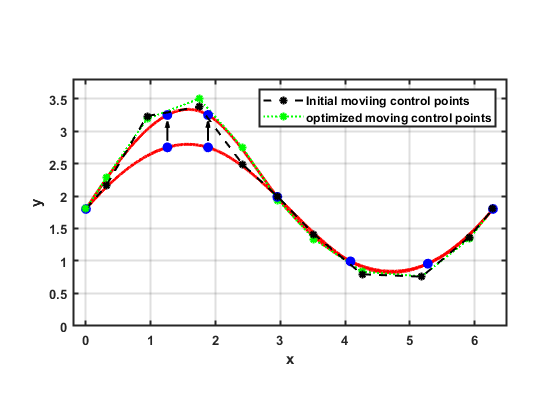}
\caption{Evolution of points on a planar curve over a time step \( \triangle t \) and control point optimization. The red curve represents the evolving geometry of the curve. The given data points (\(\vec{q}_i\)) are labeled as \(\blue{\bullet}\), with two points dynamically moved based on geometric quantities. The initial positions of the control points (\(\vec{P}_j\)) are labeled as \(*\), and their optimized positions, which minimize the interpolation error, are labeled as \(\green{*}\).}
\label{fig:movingpoints}
    \end{figure}

Once the geometric properties of the point cloud are known, these quantities directly influence the motion of the points on the evolving curve. Instead of remaining static, each point \(\vec{q}_i(t) \in \mathbb{R}^2\) moves dynamically according to the velocity field, which governs the deformation of the point cloud over time. Specifically, this motion is curvature-driven and follows the formulation introduced in Section~\ref{sec:general_framework}, where the velocity field determines the updated position of each point. The process is visually illustrated in Figure~\ref{fig:movingpoints}, which highlights the initial and updated positions of the data points and the general movement of the curve over a time step.

 Similarly, the control points  
 \(\{\vec{P}_j\}_{j=1}^m\) 
 of the curve \(\vec{f}(u)\) evolve in a manner that reflects the motion of points on the curve. Their update follows the same principle as the data point evolution, where the velocity field drives their motion over time:

\[
\vec{P}_j(t + \triangle t) = \vec{P}_j(t) + \triangle t \cdot \vec{V}_j(t).
\]
Once the core points
\(\vec{q}_i \in \mathcal{N}_{k(i)}^{core}\)
and the control points \(\{\vec{P}_j\}_{j=1}^m\) are updated, the interpolation error $err_{\text{interp}}$ is evaluated. The interpolation error is defined as the maximum distance between the points \(\vec{q}_i \in \mathcal{N}_{k}\) and the corresponding points on the curve \(\vec{f}(u)\):
\[
err_{\text{interp}} = \max_{\vec{q}_i \in N_{m}} \left\| \vec{q}_i - \vec{f}(u_i) \right\|.
\]
If \(err_{\text{interp}} > \tau\), where \(\tau\) is a predefined tolerance, we proceed to optimize the control points \(\{\vec{P}_j\}_{j=1}^m\) to minimize the interpolation error. This optimization is achieved by solving the following minimization problem:
\[
\min_{\{\vec{P}_j\}} \sum_{\vec{q}_i \in  \mathcal{N}_{k}} \left\| \vec{q}_i - \vec{f}(u_i) \right\|^2.
\]
An iterative method, such as Gauss-Seidel, is employed to update the control points. The update rule for each iteration is given by:

\[
\vec{P}_j^{(k+1)} = \vec{P}_j^{(k)} - \alpha \nabla_{\vec{P}_j} \sum_{\vec{q}_i \in  \mathcal{N}_{k}} \left\| \vec{q}_i - \vec{f}(u_i) \right\|^2,
\]
where \(k\) is the iteration index, \(\alpha\) is the step size, and \(\nabla_{\vec{P}_j}\) is the gradient of the objective function with respect to \(\vec{P}_j\).

This process ensures that both the core points \(\vec{q}_i \in \mathcal{N}_{k(i)}^{core}\) and the control points \(\{\vec{P}_j\}_{j=1}^m\) evolve in a manner that balances the smoothness of the curve with the accuracy of the interpolation. As a result, the curve remains aligned with the updated point cloud, providing a precise geometric representation of the evolving data.  Each time step \( \triangle t \) updates the configuration of points \( \{ \vec{q}_i(t + \triangle t) \}_{i=1}^N \), driving the evolution of the manifold. As shown in Figure~\ref{fig:movingpoints}, the red curve represents the evolving geometry of the manifold. The given data points \(\vec{q}_i\), labeled as \blue{$ \bullet$}, are moved according to the geometric quantities. The updated positions of the initial control points \(\vec{P}_j\), labeled as \(*\), also following the geometric quantities. The interpolation error is evaluated after the points move, and an iterative method is applied to optimize the control points, as indicated by \(\green{*}\). This framework is particularly effective for surface reconstruction, real-time simulations, and modeling deformable shapes \cite{olshanskii2017trace}.
\subsection{Point removal, insertion, and redistribution} \label{sec:rem_ins_distr}

After updating the points \(\vec{q}_i(t)\) and their control points \(\vec{P}_j(t)\), it is crucial to maintain an appropriate spacing of the points along the B-Spline curve \(\vec{f}(u)\). Balanced spacing ensures accurate interpolation while avoiding over-sampling or undersampling in specific regions. This process involves three key steps: removing points that are excessively close, inserting points where spacing is insufficient, and redistributing points to achieve near-uniform spacing.

Let \(d_i\) denote the Euclidean distance between two consecutive points \(\vec{q}_i\) and \(\vec{q}_{i+1}\), as previously defined. If this distance satisfies \(d_i < d_{\text{tol,min}}\), one of the points (typically \(\vec{q}_{i+1}\)) is removed from the point cloud to prevent oversampling in that region. Conversely, if \(d_i > d_{\text{tol,max}}\), a new point is inserted between \(\vec{q}_i\) and \(\vec{q}_{i+1}\). The new point is placed along the curve \(\vec{f}(u)\) by selecting an intermediate parameter value \(u_{\text{new}}\) such that the corresponding point \(\vec{q}_{\text{new}} = \vec{f}(u_{\text{new}})\) reduces the gap between \(\vec{q}_i\) and \(\vec{q}_{i+1}\). These adjustments ensure that no two consecutive points are too close or too far apart.

Once the point cloud has been refined through removal and insertion, the points are redistributed along the B-Spline curve \(\vec{f}(u)\) to achieve a near-uniform spacing. The total arc length \(L\) of the curve is approximated numerically as:

\[
L \approx \sum_{i=1}^{m-1} \|\vec{f}(u_{i+1}) - \vec{f}(u_i)\|,
\]
where \(u_i\) are the parameter values corresponding to the points \(\vec{q}_i\). Based on the total arc length \(L\) and the updated number of points \(m\), the target spacing between consecutive points is calculated by $d_{\text{target}} = \frac{L}{m}$.
To achieve this target spacing, the values of the new parameters \(u_i'\) are determined such that the redistributed points \(\vec{q'}_i = \vec{f}(u_i')\) satisfy the following:

\[
\left| \|\vec{q'}_{i+1} - \vec{q'}_{i}\| - d_{\text{target}} \right| < \varepsilon_d,
\]
where \(\varepsilon_d\) is a tolerance that controls the allowed deviation from the target spacing. This redistribution ensures that the points are nearly equidistant along the curve.

The removal, insertion, and redistribution process of points is applied iteratively to all points \(\vec{q}_i \in Q\), where each point belongs to the stencil \(\mathcal{N}_{k}\) as defined in \eqref{eq:cover}. After each iteration, the convergence criteria are evaluated. The algorithm continues until the point cloud satisfies the desired configuration or the final time has been reached.

\subsection{Computational cost analysis of evolving point cloud data}\label{sec:cosanalysis}
The computational cost of frequent updates to point positions and geometric properties is a significant challenge in dynamic simulations of evolving point clouds sampled from smooth manifolds. This section compares the traditional point evolution method with the proposed approach, which improves computational efficiency by utilizing local B-spline interpolation based on control points.

In any non-updated version where one stencil is used for each data point, the cost of estimating geometric quantities between neighboring data points is determined by the cost of interpolation. At \( t = 0 \), the initial setup involves computing control points for the \( N \) data points. For each point, interpolation is performed based on its local neighborhood, which includes \( m \) neighboring points. This process incurs a computational complexity of \( \mathcal{O}(N  m^3) \), where \( m \) represents the stencil size around each point. As the point cloud evolves over the time interval \( [0, T] \), the traditional method recalculates all control points at each time step \( \Delta t \). Consequently, the cumulative cost of this re-interpolation process is \( \mathcal{O}(N m^3 T/\Delta t) \), making the method computationally expensive and inefficient for high-resolution simulations requiring small \( \Delta t \).

Using one stencil for multiple data points reduces computational costs by optimizing the computation of control points. Calculating the initial control points incurs an interpolation cost of \( \mathcal{O}(m^3) \) per stencil. Partitioning the point cloud into approximately \( N / m \) non-overlapping \emph{core covers} reduces the total computational cost for the initial setup, reducing it from \( \mathcal{O}(N m^3) \) to \( \mathcal{O}(N m^2) \). During evolution (\( t + \Delta t \)), instead of re-computing all control points, only the affected ones within the relevant stencil \( \mathcal{N}_m \) are updated iteratively. This adjustment reduces the cost per time step from \( \mathcal{O}(N m^3) \) to the more efficient \( \mathcal{O}(N m^2) \). During the time interval \( [0, T] \), the cumulative cost of the updates becomes \( \mathcal{O}(N m^2 T / \Delta t) \), offering a significant improvement over the traditional cumulative cost of \( \mathcal{O}(N m^3 T / \Delta t) \). Combining the initial setup and iterative updates, the total computational cost of the proposed approach is \( \mathcal{O}(N m^3 + N m^2 T / \Delta t) \), demonstrating substantial efficiency gains compared to the traditional method. Actual runtime may vary depending on implementation and computational environment.

\section{Numerical Results and Simulations}\label{sec:results}
In this section, we present and analyze the numerical results obtained by simulating the evolution of the curve driven by the velocity field \( \vec{V} \). The simulation relies on B-Spline interpolation to approximate the geometric quantities, with several key parameters that affect accuracy and stability: the degree $p$ of the B-Spline basis functions, the stencil size $m$, and the density of the point cloud \( h \).

The temporal resolution, controlled by the time step $\Delta t$, plays a critical role in the stability and accuracy of simulations. As in many numerical methods, finding an appropriate time step is a common challenge, as it directly affects both stability and computational efficiency. In all examples presented here, we employ a fixed time step $\Delta t$ chosen to be sufficiently small to ensure numerical stability and accuracy. We note that the selection of the time step, as well as the choice of time integration scheme (explicit, implicit, or IMEX), is independent of our geometric evolution framework and can be adapted for specific applications. A systematic investigation of optimal or adaptive time stepping strategies, though important, is outside the scope of this work.

To evaluate the algorithm performance, we systematically examine how varying the B-Spline degree, stencil size, and point cloud density  influences the interpolation error   $err_{\text{interp}}$, which is compared with the predefined tolerance $\tau$
to ensure reliable approximations of the geometric quantities and maintain numerical stability. The results demonstrate how these parameters affect both simulation accuracy and computational efficiency, while also highlighting the trade-offs between accuracy and cost, providing practical insights into selecting optimal configurations for different simulation scenarios.
\begin{figure}[t!]
\begin{center}
\subfigure[]{\includegraphics[width=2.9in]{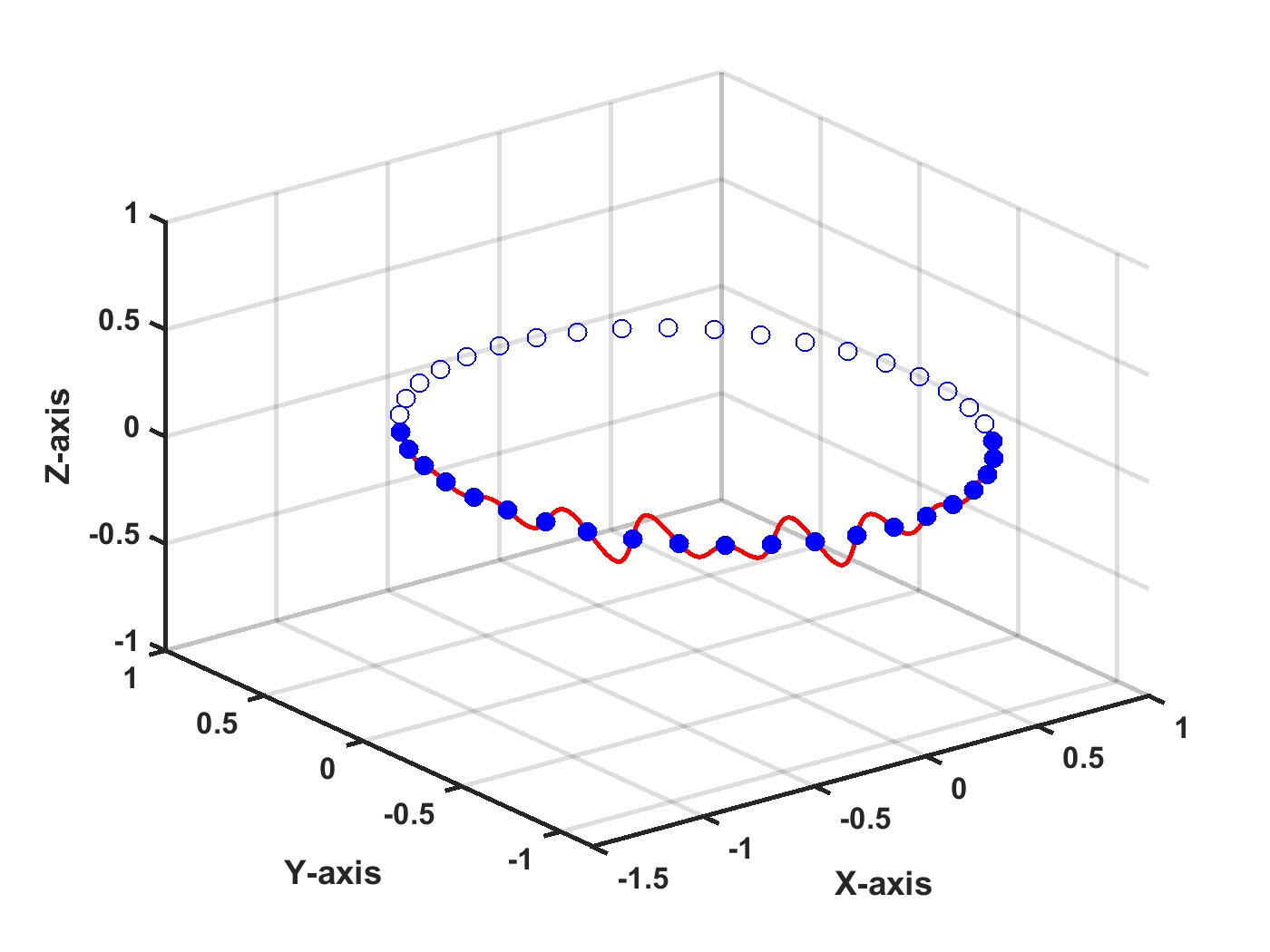}}
 \subfigure[]{\includegraphics[width=2.9in]{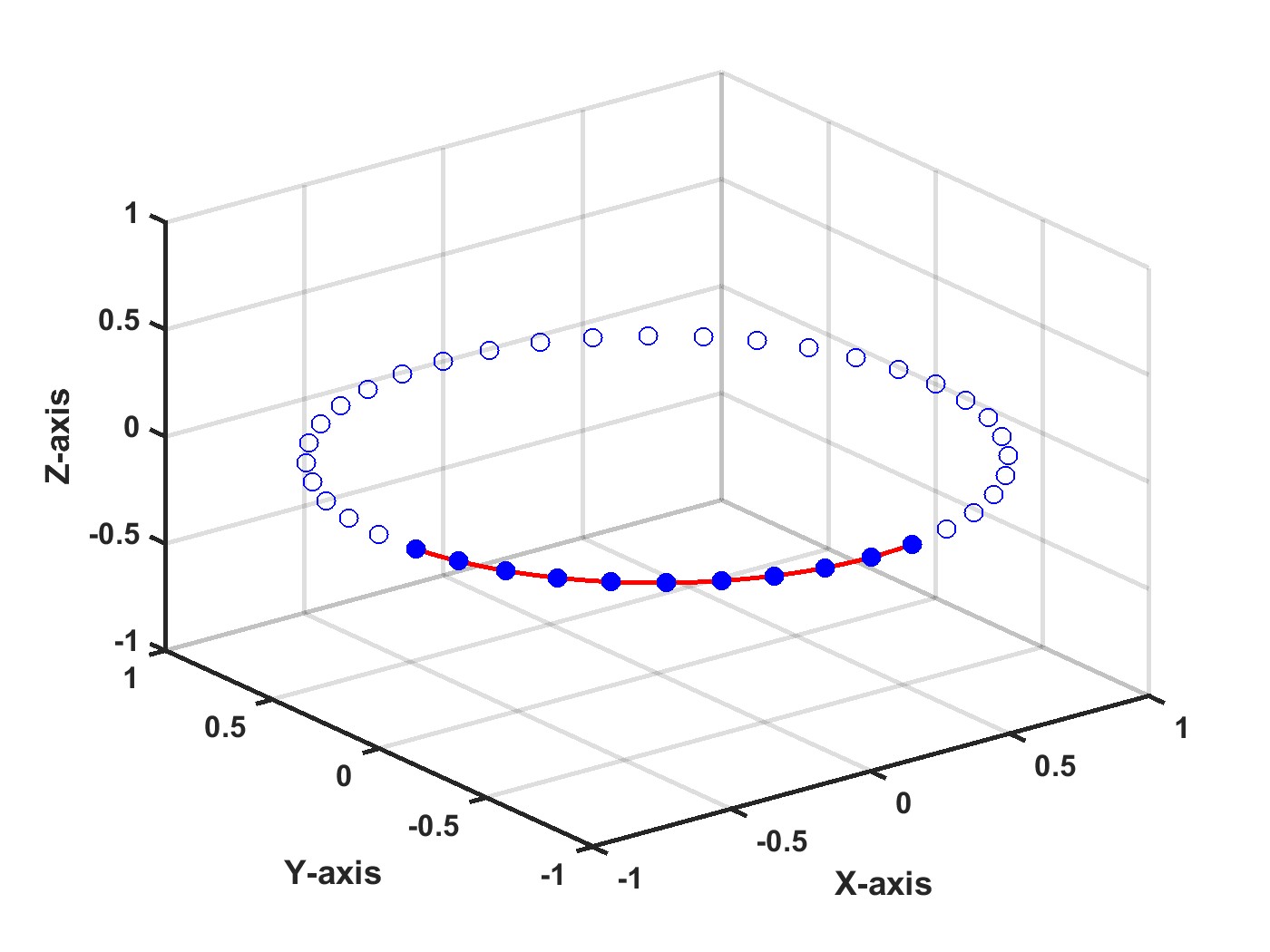}}
\caption{Impact of degree and stencil size  on B-Spline interpolation (a) For a large stencil size: when $m=20$ at degree $p= 3 $
(b) For a smaller stencil size: when $m=10$ at $p = 3.$ }
\label{fig:1d_interpolation}
\end{center}
\end{figure}

\begin{figure}[ht!]
\begin{center}
\subfigure[Normal error variation with B-Spline degree and stencil size ]{\includegraphics[width=3in]{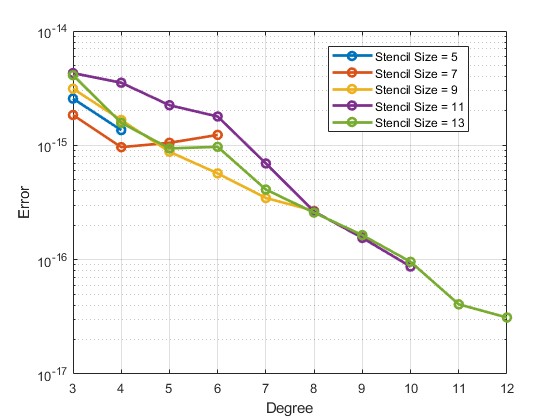}}
\subfigure[Curvature error variation with B-Spline degree and stencil size]{\includegraphics[width=3in]{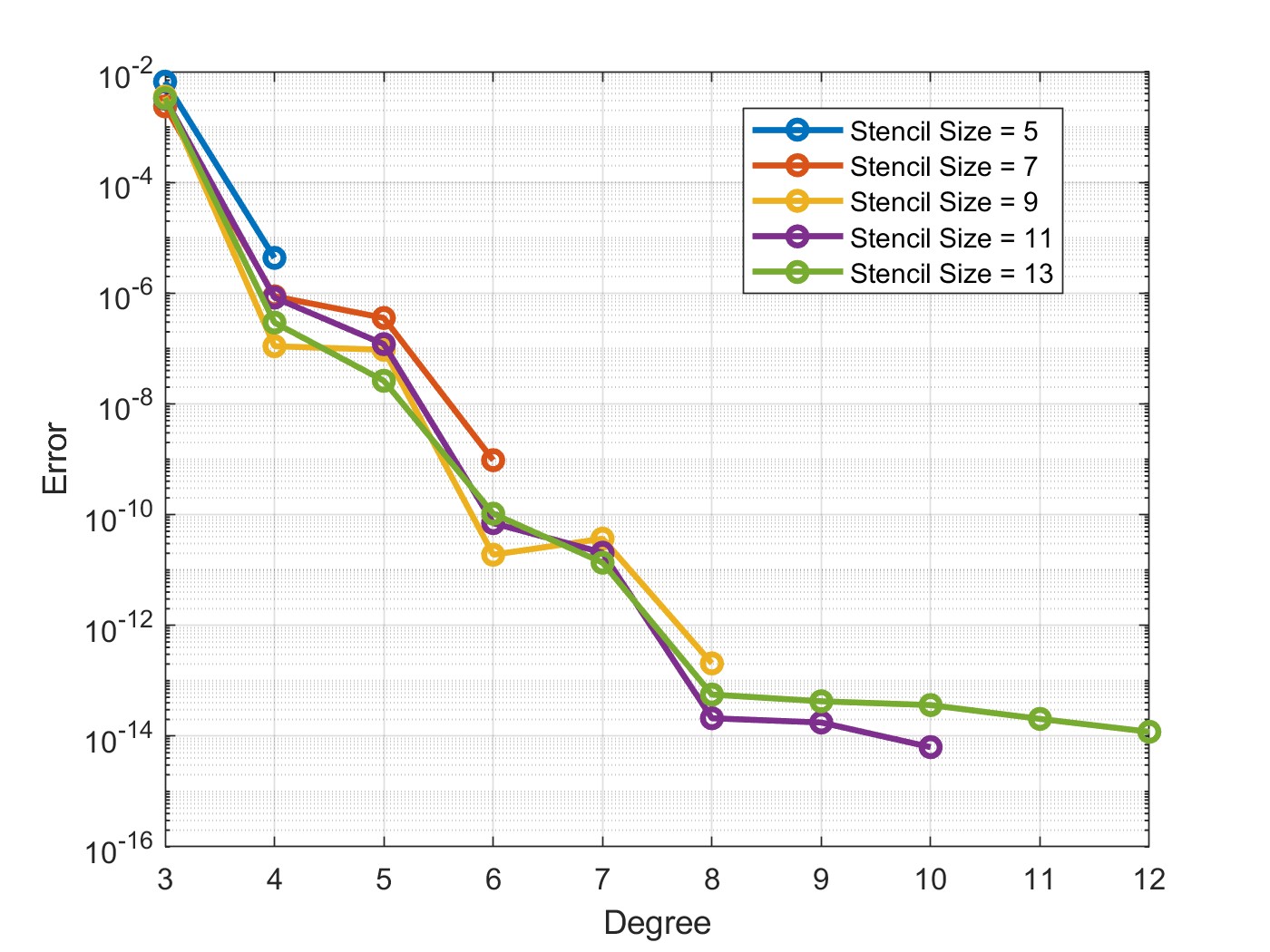}}
\caption{ Optimal parameter selection for algorithm tuning: (a) Normal error as a function of B-Spline degree and stencil size, and (b) Curvature error as a function of B-Spline degree and stencil size. For B-Splines, the maximum allowable degree is generally limited to \(m - 1\), where \(m\) is the stencil size. As a result, the error curves are shorter for smaller stencil sizes because higher B-Spline degrees are not feasible. Conversely, larger stencil sizes allow for higher degrees, resulting in longer error curves.}
\label{fig:dvsstensil}
\end{center}
\end{figure}

\subsection*{Example 1: Impact of stencil size, B-Spline degree, and point cloud density on control point evolution}
In this example, we explore how the choice of stencil size, B-Spline degree, and point cloud density affects the evolution of control points during the interpolation process. Consider a set of $ N = 40 $ points sampled from a circle. B-spline interpolation is carried out with varying stencil sizes and degrees to assess their influence on the smoothness of the curve, as well as the precision of the normal $ \n $, curvature $ \kappa $, and the behavior of the control points.
For a large stencil size, $ m = N/2 = 20 $, the interpolated curve exhibits significant oscillations, as seen in Figure \ref{fig:1d_interpolation}(a), leading to instability and inaccurate normal and curvature estimations. This causes the control points $\vec{P}_j$ to misalign with the data points $\vec{q}_j$, making them difficult to manage. In contrast, with a smaller stencil size, $ m = 10 $, the curve becomes smoother, as shown in Figure \ref{fig:1d_interpolation}(b), with reduced oscillations and improved normal and curvature accuracy. As a result, the control points stay aligned with the data points, offering better control over the curve's shape and ensuring a more stable interpolation.

To further investigate the impact of various stencil sizes and degrees on interpolation accuracy, we extend the analysis in Figure \ref{fig:dvsstensil}, where the errors in the normal and curvature estimations are plotted for different degrees and stencil sizes. This broader examination reveals that a mid-range stencil size, such as $ m = 9 $, coupled with degrees $ m-1 $ or $ m-2 $, consistently minimizes the errors. This combination not only improves the stability of the curve, but also ensures smoother and more predictable control point evolution, making it well-suited for precise interpolation tasks. Larger stencil sizes introduce instability, while smaller ones may lack the detail needed to accurately capture the curve's complexity, highlighting the importance of selecting an appropriate balance between stencil size and degree for optimal interpolation results.

\begin{figure}[t!]
\begin{center}
\subfigure[Normal error variation with stencil size and point cloud density]{\includegraphics[width=3.1in]{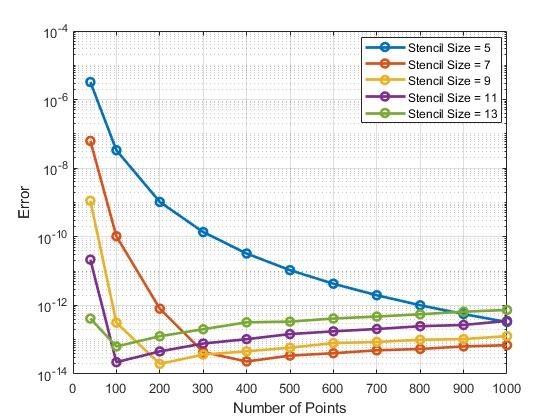}}
\subfigure[Curvature error variation with stencil size and point cloud density]{\includegraphics[width=3.1in]{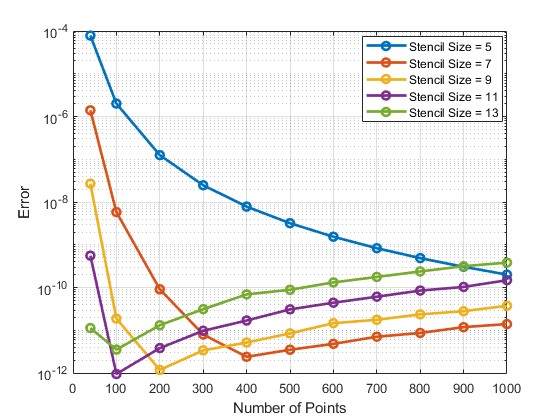}}
\caption{Optimal parameter selection for algorithm tuning: (a) Variation of normal error as a function of stencil size and point cloud density, and (b) Variation of curvature error as a function of stencil size and point cloud density.}
\label{fig:stensilvsdesnity}
\end{center}
\end{figure}

As the number of points increases to a denser dataset, up to $ N = 1000 $, as shown in Figure \ref{fig:stensilvsdesnity}, the choice of stencil size remains pivotal. Despite the higher point density, the mid-range stencil size $ m = 9 $ provides the best balance, minimizing normal and curvature errors. The result illustrates the efficiency and adaptability of the mid-range stencil size, which performs well with smaller datasets and seamlessly accommodates denser point clouds. These results show that careful selection of stencil size and degree is key to achieving stable and accurate interpolation.
\subsection*{Example 2: Simulation of interface dynamics with velocity $\vec V =-\kappa \n $}
Consider a circular interface with an initial radius \( r_0 = 1 \). Under the velocity field \( \vec{V} = -\kappa \n \), the interface evolves over time and the radius decreases. Our method starts with $N = 200$ discretized points along the interface. As the radius decreases, the number of points is dynamically reduced to avoid overcrowding, particularly in low-curvature regions. For instance, at time $t = 0.2545$, the number of points is reduced to $N = 99$, and by time $t = 0.4495$, the number of points is further reduced to $N = 49$. This reduction in points is efficiently handled by the proposed B-Spline interpolation method, which automatically manages point removal and redistribution throughout the process. The method ensures smoothness and accuracy despite the reduction in points, significantly reducing computational costs while maintaining geometric fidelity. The technique excels at efficiently redistributing points to maintain uniform spacing as the circle shrinks and the curvature evolves, making it particularly well-suited for curvature-driven flows.

\begin{figure}
    \centering
  \subfigure[] {\includegraphics[width=0.45\textwidth]{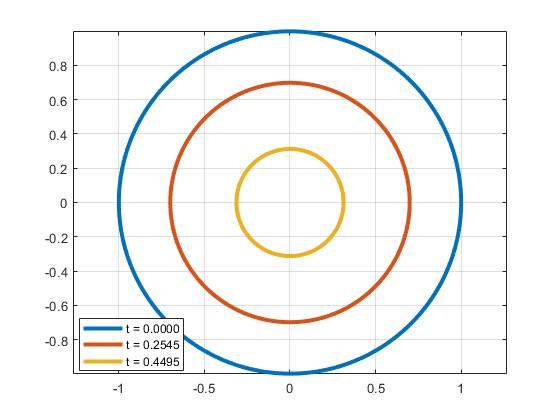}}
  \subfigure[] {\includegraphics[width=0.45\textwidth]{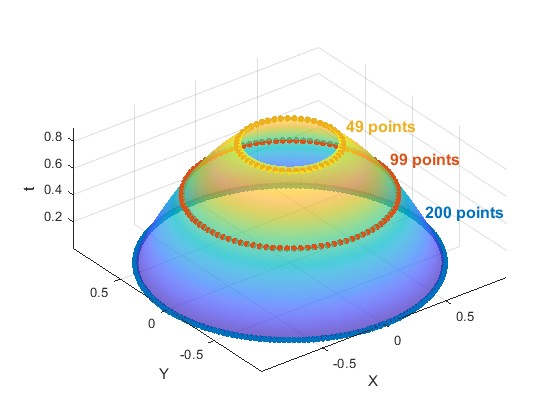}}
    \caption{Evolution of a circular interface under curvature-driven flow. (a) 2D view at $t = 0.0$, $t = 0.2545$, and $t = 0.4495$. (b) 3D view illustrating the interface shrinking over time, with the number of points reducing from 200 to 99, and finally to 49, while maintaining smoothness and accuracy.}
    \label{fig:circle}
\end{figure}

An example result of the shrinking circular interface is shown in Figure \ref{fig:circle}, illustrating the evolution of the interface at different time steps.
The collapse of a circular interface can be described analytically, making it an ideal test case to assess the accuracy of our method. The radius of a circle at time $t$ is given by the equation:
\begin{equation}
r(t) = \sqrt{r_0^2 - 2t},  \quad t \in [0, 2\pi].
\label{eq:circle_test}
\end{equation}
where $ r_0$ is the initial radius. The evolution of the circle's radius was tracked until $ t = 0.180 $, using varying arc length spacings $ h $ and a time step $ \triangle t = 1 \times 10^{-3} $. The accuracy of the B-Spline and PHS-RBF+Poly interpolation methods in approximating the radius was compared to the exact value, and the resulting errors are summarized in Table \ref{tab:comparison}. Both methods demonstrate comparable accuracy across different values of $ N $ and $ h $, achieving similar $ L_2 $ errors, especially as $ N $ increases and $ h $ decreases. However, a key distinction lies in the efficiency of the B-Spline method, which attains this accuracy with a smaller stencil size. The B-Spline method uses stencil sizes of $ m = 5 $ or $ m = 9 $, depending on configuration. In contrast, the PHS-RBF+Poly method, which incorporates a radial basis function of order 3 combined with a polynomial of degree 3, requires a larger stencil. Results for this method are shown for $m = 15$. 

To further clarify the comparative performance of these methods, we systematically examined both B-Spline and PHS-RBF+Poly approaches over a wide range of time step sizes, from $1 \times 10^{-5}$ up to $0.5$. In the B-Spline framework, the control points are updated adaptively; at each time step, the evolved control points provide near-optimal initial guesses for the data points. This efficient initialization enables the Gauss-Seidel method to converge in a single iteration for all time steps up to $\Delta t = 0.1$, and only two iterations for $\Delta t = 0.5$. For the PHS-RBF+Poly method in our test configuration, the previous time step's coefficients do not provide as effective an initial guess, and iterative solvers such as conjugate gradient require additional iterations to achieve convergence.

A detailed analysis of the interpolation matrices reveals distinct characteristics of each approach. The B-Spline interpolation matrix exhibits a condition number of approximately $11.8$, a maximum eigenvalue of one, and smoothly scaling changes in the control points ($|\Delta P|_2$) with values ranging from $10^{-6}$ to $10^{-1}$ across all tested time steps. The PHS-RBF+Poly interpolation matrix in our configuration shows different properties, with a higher condition number and larger spectral radius. While both methods can be effectively solved using appropriate iterative techniques, the key advantage of the B-Spline method lies in its ability to use evolved control points as high-quality initial guesses for subsequent time steps, enabling rapid convergence even with compact stencils.

\newcolumntype{R}[1]{>{\raggedleft\arraybackslash}p{#1}}
\newcolumntype{C}[1]{>{\centering\arraybackslash}p{#1}}
\begin{table}[t!]
\caption{Comparison of $L_2$-error in circle radius accuracy for B-Spline and PHS-RBF methods under the velocity field $\vec{V} = -\kappa \n$, evaluated for different values of $N$ and $h$.}
\centering
\begin{tabular}{c c c c c}
\toprule
$N$ & $h$ & B-Spline ( $m = 5$) & B-Spline ( $m = 9$) & PHS-RBF+Poly ($m = 15$) \\
\midrule
30  & $2.09\times 10^{-1}$ & $9.11\times 10^{-4}$ & $8.39\times 10^{-4}$ & $8.29\times 10^{-4}$  \\
60  & $1.05\times 10^{-1}$ & $4.84\times 10^{-4}$ & $4.35\times 10^{-4}$ & $4.17\times 10^{-4}$ \\
120 & $5.24\times 10^{-2}$ & $2.46\times 10^{-4}$ & $2.23\times 10^{-4}$ & $2.09\times 10^{-4}$  \\
240 & $2.62\times 10^{-2}$ & $1.60\times 10^{-4}$ & $1.35\times 10^{-4}$ & $1.04\times 10^{-4}$  \\
\bottomrule
\end{tabular}
\label{tab:comparison}
\end{table}

For more complex geometries, the dynamics of curvature-driven flow become even more pronounced. A prime example of this is a two-dimensional asterisk-shaped interface, which exhibits significantly sharper features and higher variations in curvature compared to simpler shapes like circles.

\begin{figure}[ht!]
    \centering
  \subfigure[] {\includegraphics[width=0.45\textwidth]{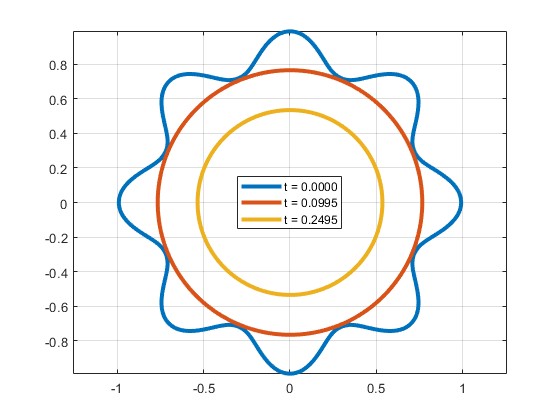}}
  \subfigure[] {\includegraphics[width=0.45\textwidth]{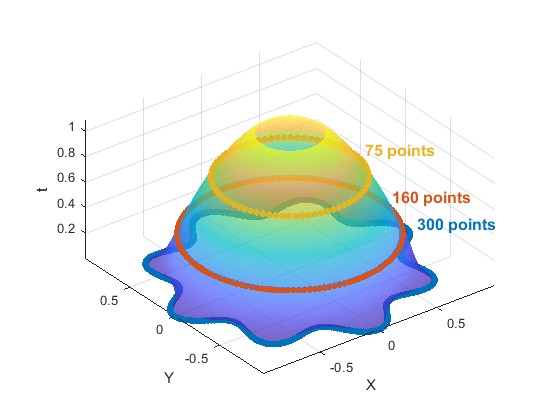}}
    \caption{Evolution of a two-dimensional asterisk-shaped interface under curvature-driven flow at different time steps. (a) 2D view of the asterisk's evolution. (b) 3D representation of the asterisk's evolution, where the height corresponds to a scaled version of time for better visualization.}
    \label{fig:asterisk}
\end{figure}

In this case, the radius of the asterisk-shaped interface is given by
\begin{equation}
r = 1 + 0.3 \cos^2(4t),  \quad t \in [0, 2\pi].
\end{equation}
and its evolution is depicted in Figure \ref{fig:asterisk} at various time steps: $t = 0,0.0995$, and $ 0.2495$. At $t = 0$, the interface starts with a distinct asterisk-like configuration, characterized by sharp, pronounced points. As time progresses, the curvature-driven flow leads to a smoothing of these points. The regions with higher curvature evolve more rapidly, causing sharp features to gradually diminish. At later times, the shape becomes increasingly circular as the curvature effects drive the interface toward uniformity.

A key strength of the proposed B-Spline interpolation method is its ability to dynamically reduce the number of points as the shape evolves, particularly in regions of low curvature, while maintaining smoothness and preserving geometric integrity. As shown in the 3D representation (Figure \ref{fig:asterisk}b), the number of points decreases from 300 initially to 160 at \( t = 0.0995 \), and to 75 at \( t = 0.2495 \), all while accurately representing the evolving interface.  This automatic point reduction is essential in curvature-driven flows, allowing the algorithm to adapt the point distribution for capturing rapid changes in high-curvature regions while maintaining accuracy in smoother areas. The asterisk interface exemplifies the method's ability to balance computational efficiency and geometric fidelity, even for complex shapes.
\subsection*{Example 3: Evolving circle with coupled reaction-diffusion}

In this example, we investigate the evolution of a circular boundary, which changes over time on the basis of the solution of a reaction-diffusion system defined on the circle. The system is modeled by two coupled PDEs on the evolving boundary, denoted by $ \Omega(t) $, representing the circle. The governing equations are:
\begin{align}
    \frac{\partial u}{\partial t} &= D_u \Delta_{\Omega} u - u (\nabla_{\Omega} \cdot \vec V) + g_1(u, v), \\
    \frac{\partial v}{\partial t} &= D_v \triangle_{\Omega} v - v (\nabla_{\Omega} \cdot \vec V) + g_2(u, v),
\end{align}
where $ u $ and $ v $ are scalar fields defined on the evolving circle $ \Omega(t) $. The operator $ \triangle_{\Omega} $ denotes the Laplace-Beltrami operator on the circle, and $ \vec V $ is the velocity field that governs the time evolution of the boundary. The terms $ g_1(u, v) $ and $ g_2(u, v) $ represent the reaction kinetics, while $ D_u $ and $ D_v $ are the diffusion coefficients for $ u $ and $ v $.

The velocity $ \vec V $ of the evolving boundary is now given by:
\begin{equation}
    \vec V = \left(c_1 \kappa + c_2 u\right) \n,
\end{equation}
where $ \kappa $ is the curvature of the evolving boundary, and $ \n $ is the outward-pointing unit normal vector. The velocity depends on the solution of the PDE $u$, making the system strongly coupled, since the evolution of the boundary and the dynamics of reaction-diffusion are interdependent. The parameters $ c_1 $ and $ c_2 $ control the relative contributions of the curvature-driven flow and the reaction-diffusion field to the boundary evolution.
\begin{figure}[t!]
    \centering
  \subfigure[] {\includegraphics[width=0.42\textwidth]{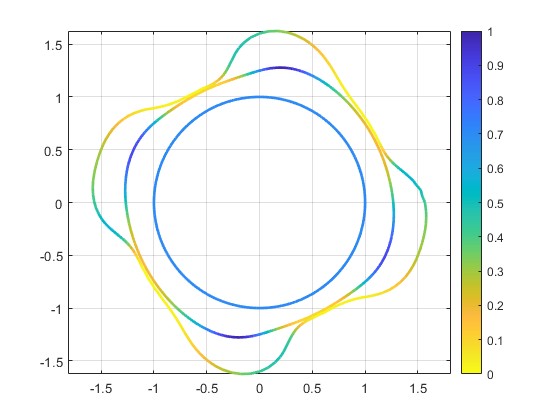}}
 \subfigure[] { \includegraphics[width=0.42\textwidth]{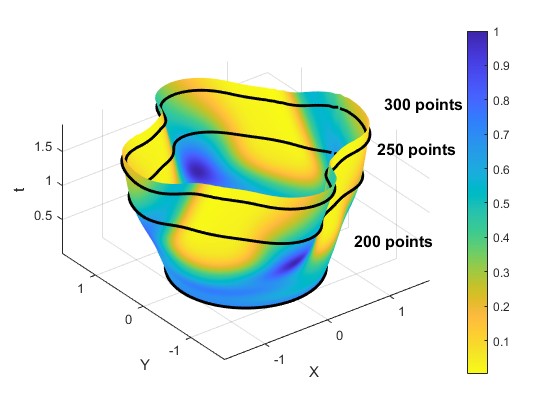}}\\
  \subfigure[] {\includegraphics[width=0.42\textwidth]{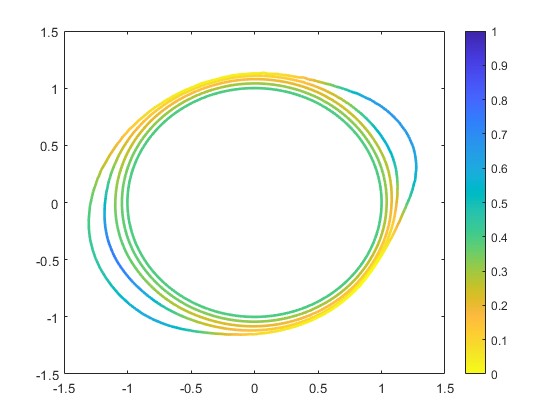}}
 \subfigure[] { \includegraphics[width=0.42\textwidth]{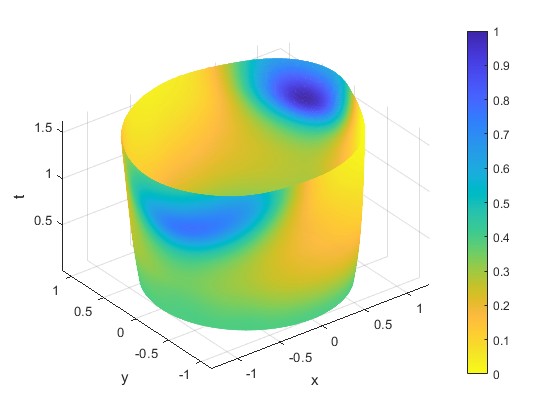}}\\
    \caption{Effect of varying diffusion coefficients on pattern evolution. Left: 2D views of the evolving boundary, where the color indicates the concentration of $u(\theta, t)$. Right: 3D surface representations where the height corresponds to the concentration of $u(\theta, t)$. (Top) For $D_u=0.1$ and $D_v=1.5$, smoother patterns and more pronounced boundary deformations are observed. (Bottom) For $D_u=0.1$ and $D_v=0.4$, sharper patterns with less boundary deformation are shown.}
    \label{fig:pattern_evolution}
\end{figure}

For the reaction terms, we set:
\begin{align*}
    g_1(u, v) &= \gamma (c - u + u^2 v), \\
    g_2(u, v) &= \gamma (d - u^2 v),
\end{align*}
where $ \gamma $ is a reaction rate constant, and $ c $ and $ d $ are parameters that control the reaction kinetics. The initial condition is chosen as a small Gaussian perturbation added to the homogeneous steady state. Specifically, we define the initial conditions for $ u $ and $ v $ as:
\begin{align*}
    u(\theta, 0) &= u_0 + 0.5 u_0 \exp\left(-\frac{(\theta - \theta_0)^2}{2 \sigma^2}\right), \\
    v(\theta, 0) &= v_0 + 0.5 v_0 \exp\left(-\frac{(\theta - \theta_0)^2}{2 \sigma^2}\right),
\end{align*}
where $ \theta $ is the angular position on the circle, $ u_0 = c + d $, and $ v_0 = \frac{d}{(c + d)^2} $ are the homogeneous steady-state values. The parameter $ \sigma $ controls the width of the Gaussian perturbation, and the perturbation is centered at a random angle $ \theta_0 $, introducing randomness to facilitate the formation of non-uniform patterns.

We discretize the circle into 200 uniformly distributed points and compute the Laplacian \( \triangle_{\Omega} \) using a multiscale radial basis function (RBF) method. Time-stepping is performed using an implicit-explicit (IMEX) scheme \cite{ascher1995implicit}, where diffusion terms are treated implicitly, and nonlinear reaction terms \( g_1(u, v) \) and \( g_2(u, v) \) explicitly. Figure \ref{fig:pattern_evolution} illustrates the evolution of boundary shapes and reaction-diffusion patterns for varying diffusion coefficients \( D_v \), demonstrating the algorithm’s capacity to adjust the point density during deformation dynamically.
\begin{figure}[t]
    \centering
  \subfigure[$c_1 = 0.02, c_2 = 1$] {\includegraphics[width=0.33\textwidth]{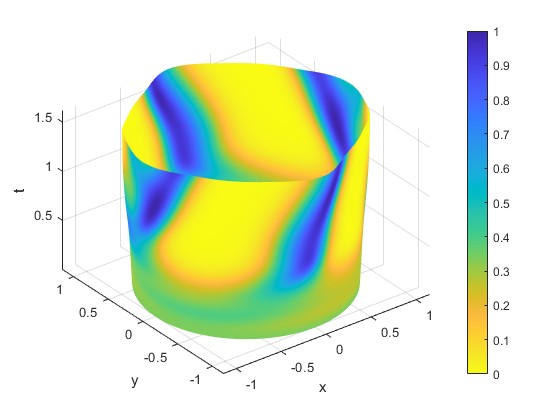}}
 \subfigure[$c_1 = 0.03, c_2 = 5$] { \includegraphics[width=0.32\textwidth]{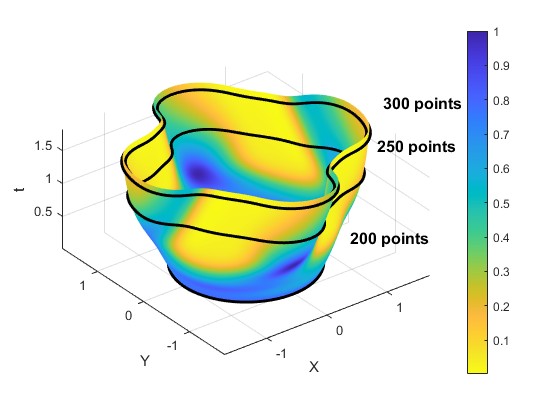}}
  \subfigure[$c_1 = 0.03, c_2 = 10$] {\includegraphics[width=0.32\textwidth]{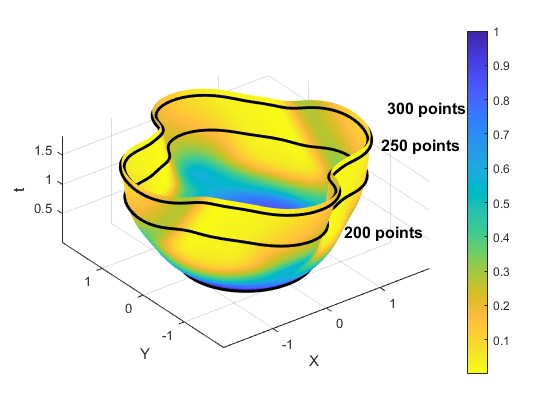}}
    \caption{Effect of the solution $u$ on boundary velocity $\vec V$. As the coupling parameter $c_2$ increases, the influence of the reaction-diffusion solution $u$ on boundary velocity becomes more significant. Higher values of $c_2$ result in sharper patterns and more complex boundary deformations due to the stronger contribution of $u$ in the velocity term $\vec V = (c_1 \kappa + c_2 u) \n$.}
    \label{fig:pattern_evolution2}
\end{figure}
In Figures \ref{fig:pattern_evolution}(a) and \ref{fig:pattern_evolution}(c), two-dimensional projections show the boundary evolution for \( D_v = 1.5 \) and \( D_v = 0.4 \), respectively. For \( D_v = 1.5 \), the boundary undergoes significant deformation, forming irregular shapes with sharp features corresponding to the localized reaction-diffusion patterns in \( u \). The algorithm increases the point count from 200 to 300 to capture these complexities. For \( D_v = 0.4 \), slower diffusion results in smoother patterns, the boundary maintaining a regular shape and the initial point density. Figures \ref{fig:pattern_evolution}(b) and \ref{fig:pattern_evolution}(d) provide three-dimensional views of these results, confirming the algorithm’s ability to balance computational efficiency with geometric accuracy by dynamically adapting the distribution of points.

Figure \ref{fig:pattern_evolution2} examines the effect of varying \( c_2 \) on boundary deformation, with \( c_1 = 0.03 \), \( D_u = 0.1 \), and \( D_v = 1.5 \) fixed. For \( c_2 = 1 \), the coupling between \( u \) and the boundary is weak, resulting in minimal deformation and smooth patterns. As \( c_2 \) increases to 10, stronger coupling intensifies boundary undulations and creates sharper features. At \( c_2 = 10 \), the reaction-diffusion dynamics dominate, producing irregular boundary shapes with sharp gradients and significant deformation. These results illustrate that increasing \( c_2 \) amplifies boundary deformation and enhances pattern localization.

\section{Conclusion}\label{sec:conclusion}
This preliminary study examines the use of adaptive B-Spline interpolation within a Lagrangian framework to model the evolution of point cloud data. The flexibility of B-Spline basis functions allows for the efficient computation of geometric quantities, such as normal vectors and curvature, while enabling dynamic point addition and removal. By avoiding re-interpolation at each time step, the framework offers a computationally efficient alternative to traditional Radial Basis Function (RBF) techniques, particularly for problems involving evolving point clouds with varying densities. The ultimate goal of this work is to extend the approach to higher-dimensional manifolds and surfaces, with this study serving as a foundational step in that direction.

Numerical experiments validate the effectiveness of the framework through simulations of curvature-driven flows and coupled reaction-diffusion systems. For curvature-driven flows, the method accurately captured the collapse of simple geometries, such as circles, and the smoothing of more complex shapes, such as asterisks, aligning with theoretical predictions. For the coupled reaction-diffusion system, the framework successfully modeled the interplay between boundary deformation and scalar field evolution, demonstrating its ability to capture complex geometric changes and non-uniform boundary variations driven by physical processes.

While the present framework is effective for evolving closed, smooth planar curves, it assumes that the initial point cloud always represents a single, smooth manifold and does not address topological changes such as merging or splitting. However, after $t>0$, the oriented point cloud provides normal information at each point, which could potentially be used to detect when two parts of the curve approach each other from opposite sides, thus enabling the handling of such topological changes in future work (see, e.g., \cite{petras2016pdes,olshanskii2017trace}). Our experiments focus on closed curves with relatively uniform point distributions; the performance on open curves, highly non-uniform or noisy data, or curves with sharp corners and singularities remains unexplored, but these cases could be addressed by incorporating advanced point cloud processing, noise filtering, or local feature detection strategies. Future work will address these aspects, extend the approach to three-dimensional surface evolution, integrate higher-order geometric flows, and explore multi-physics applications. In particular, coupling this methodology with the Arbitrary Lagrangian-Eulerian (ALE) method~\cite{de2021numerical, mokbel2020ale} may provide further advantages for grid consistency during complex interactions such as those arising in fluid-structure dynamics.

\section*{Acknowledgment}
This work was supported by the General Research Fund (GRF No. 12301824, 12300922) of Hong Kong Research Grant Council.

  \bibliographystyle{elsarticle-num}
  \bibliography{ref-abrv}
\end{document}